\documentclass[]{article}

\newtheorem{thm}{Theorem}[section]
\newtheorem{cor}[thm]{Corollary}

\newtheorem{lem}[thm]{Lemma}
\newtheorem{Def}[thm]{Definition}
\newtheorem{rem}[thm]{Remark}

\newcommand{\be}{\begin{equation}}
\newcommand{\ee}{\end{equation}}
\newcommand{\ben}{\begin{enumerate}}
\newcommand{\een}{\end{enumerate}}

\newcommand{\qed}{\hspace*{\fill}Q.E.D.}  
\title{Homogeneous geodesics in homogeneous Finsler
spaces}
\author{%
  Dariush Latifi  \\
\\
  Faculty of Mathematics and Computer Science,\\
  Amirkabir University of Technology,\\
  424 Hafez Ave., 15914 Tehran, Iran\\
  (e-mails: dlatifi@aut.ac.ir , dlatifi@gmail.com )}

\date{}

\begin{document}
\maketitle

\begin{abstract}
In this paper, we study homogeneous geodesics in homogeneous Finsler
spaces. We first give a simple criterion that characterizes geodesic
vectors. We show that the geodesics on a Lie group, relative  to a
bi-invariant Finsler metric, are the cosets of the one-parameter
subgroups. The existence of infinitely many homogeneous geodesics on
compact semi-simple Lie group is established. We introduce the
notion of naturally reductive homogeneous Finsler space. As a
special case, we study homogeneous geodesics in homogeneous Randers
spaces. Finally, we study some curvature properties of homogeneous
geodesics. In particular, we prove that the \textbf{S}-curvature
vanishes along the homogeneous geodesics.
\end{abstract}
\\
\section{Introduction}

\

A connected Riemannian manifold $(M,g)$ is said to be
\emph{homogeneous} if a connected group of isometries $G$ acts
transitively on it. Such $M$ can be identified with
$(\frac{G}{H},g)$, where $H$ is the isotropy group at a fixed point
$o$ of $M$. The Lie algebra $\textbf{\underline{g}}$ of $G$ admits a
reductive decomposition
$\textbf{\underline{g}}=\textbf{m}\oplus\textbf{\underline{h}}$,
where $\textbf{m}\subset \textbf{\underline{g}}$ is a subspace of
$\textbf{\underline{g}}$ isomorphic to the tangent space $T_{o}M$
and $\textbf{\underline{h}}$ is the Lie algebra of $H$ [7]. In
general, such a decomposition is not unique. A \emph{homogeneous
geodesic} through the origin $o\in M=\frac{G}{H}$ is a geodesic
$\gamma(t)$ which is an orbit of a one-parameter subgroup of $G$,
that is $$\gamma(t)=exp(tZ)(o), \hskip1cm t\in R$$ where $Z$ is a
nonzero vector of $\textbf{\underline{g}}$. \

Homogeneous geodesics have important applications to mechanics. For
example, the equation of motion of many systems of classical
mechanics reduces to the geodesic equation in an appropriate
Riemannian manifold $M$.\\ Geodesics of left-invariant Riemannian
metrics on Lie groups were studied by Arnold extending Euler's
theory of rigid-body motion [1]. A major part of Arnold's paper is
devoted to the study of homogeneous geodesics. Homogeneous geodesics
are called by Arnold "relative equilibriums ". The description of
such relative equilibria is important for qualitative description of
the behaviour of the corresponding mechanical system with
symmetries. There is a big literature in mechanics devoted to the
investigation of relative equilibria. \
 Studying the set of homogeneous geodesics of a homogeneous
 Riemannian manifold $(\frac{G}{H},g)$ the concept of geodesic
 vector proved to be convenient [10]. A nonzero vector $X\in
 \textbf{\underline{g}}$ is called a geodesic vector if the curve
 $\gamma(t)=exp(tZ)(o)$ is a geodesic on $(\frac{G}{H},g)$. The
 following lemma can be found in [10].
\begin{lem}
A vector $X\in \textbf{\underline{\textsc{g}}}-\{0\}$ is a geodesic
vector if and only if $$<[X,Y]_{m},X_{m}>=0 \hskip1cm \forall
\hskip.1cm Y\in \textbf{m},$$
\end{lem}
where $< , >$ is the $Ad(H)-$invariant scalar product on
$\textbf{m}$ induced by the Riemannian scalar product on $T_{o}M$
and the subscripts $m$ indicates the projection into $\textbf{m}$.
The study of the set of homogeneous geodesics of a homogeneous
Riemannian manifold is obviously reducible to the study of the set
of its geodesic vectors.

\ A Finsler metric on a manifold is a family of Minkowski norms on
tangent spaces. There are several notions of curvature in Finsler
geometry. The flag curvature $\textbf{K}$ is an analogue of the
sectional curvature in Riemannian geometry. The Cartan torsion
$\textbf{C}$ is a primary quantity which characterizes Riemannian
metrics among Finsler metrics. There is another quantity which also
characterizes Riemannian metrics among Finsler metrics, that is the
so-called distortion $\tau$. The horizontal derivative of $\tau$
along geodesics is the so-called the \textbf{S}-curvature
$\textbf{S}= \tau_{;k}y^{k}$. While many works have been done on the
general geometric properties of Finsler geometry, such as
connections, geodesics and curvature, only very little attention has
been paid to the group aspects of this interesting field. This may
be mainly due to the reason that the Myers-Steenrod theorem in
Riemannain geometry was not successfully generalized to the
Finslerian case for a rather long period. A proof of this theorem
for the Finslerian case was given in ([19],[3]). Namely they proved
that the group of isometries of a Finsler space is a Lie
transformation group of the underlying manifold. This result opens a
door to using Lie group theory to study Finsler geometry ([4],[12]).
\

The purpose of the present paper is to study homogeneous geodesics
in homogeneous Finsler spaces. The definition of homogeneous
geodesics is similar to the Riemannian case.

\section{Preliminaries}
\subsection{Finsler spaces}

In this section, we recall briefly some known facts about Finsler
spaces. For details, see [2]. \

Let $M$ be a n-dimensional $C^{\infty}$ manifold and
$TM=\bigcup_{x\in M}
 T_{x}M$ the tangent bundle. If the continuous function $F:TM\longrightarrow
 R_{+}$ satisfies the condition that it is $C^{\infty}$ on $TM\setminus
 \{0\}$; $F(tu)=tF(u)$ for all $t\geq 0$ and $u\in TM$, i.e, $F$ is
 positively homogeneous of degree one; and for any tangent vector $y\in T_{x}M\setminus \{0\}$, the following bilinear symmetric form $g_{y}:T_{x}M \times T_{x}M\longrightarrow R$ is positive definite :
$$g_{y}(u,v)=\frac{1}{2}\frac{\partial^{2}}{\partial s\partial t}[F^{2}(x,y+su+tv)]|_{s=t=0},$$
  then we say that $(M,F)$ is a Finsler manifold.\\
Let $$g_{ij}(x,y)=(\frac{1}{2}F^{2})_{y^{i}y^{j}}(x,y).$$ By the
homogeneity of $F$, we have
$$g_{y}(u,v)=g_{ij}(x,y)u^{i}v^{j},\hskip.5cm F(x,y)=\sqrt{g_{ij}(x,y)y^{i}y^{j}}.$$

Let $\gamma:[0,r]\longrightarrow M$ be a piecewise  $C^{\infty}$
curve. Its integral length is defined as
$$L(\gamma)=\int_{0}^{r}F(\gamma(t),\dot{\gamma}(t))dt.$$ For $x_{0},x_{1}\in
M$ denote by  $\Gamma(x_{0},x_{1})$  the set of all piecewise
$C^{\infty}$ curve $\gamma:[0,r]\longrightarrow M$ such that
$\gamma(0)=x_{0}$ and $\gamma(r)=x_{1}$. Define a map $d_{F}:M\times
M\longrightarrow [0,\infty)$ by \[d_{F}(x_{0},x_{1})=\inf_{\gamma\in
\Gamma(x_{0},x_{1})} L(\gamma).\] Of course we have
$d_{F}(x_{0},x_{1})\geq 0$, where the equality holds if and only if
$x_{0}=x_{1}$; $d_{F}(x_{0},x_{2})\leq
d_{F}(x_{0},x_{1})+d_{F}(x_{1},x_{2})$. In general, since $F$ is
only a positive homogeneous function, $d_{F}(x_{0},x_{1})\neq
d_{F}(x_{1},x_{0})$, therefore $(M,d_{F})$ is only a non-reversible
metric space.

Let $\pi^{\ast}TM$ be the pull-back of the tangent bundle $TM$ by
$\pi: TM\setminus\{0\}\longrightarrow M$. Unlike the Levi-Civita
connection in Riemannian geometry, there is no unique natural
connection in the Finsler case. Among these connections on
$\pi^{\ast}TM$, we choose the \emph{Chern connection} whose
coefficients are denoted by $\Gamma^{i}_{jk}$(see[2,p.38]). This
connection is almost $g-$compatible and has no torsion. Here
$g(x,y)=g_{ij}(x,y)dx^{i}\otimes
dx^{j}=(\frac{1}{2}F^{2})_{y^{i}y^{j}}dx^{i}\otimes dx^{j}$ is the
Riemannian metric on the pulled-back bundle $\pi^{\ast}TM$.

The Chern connection defines the covariant derivative $D_{V}U$ of a
vector field $U\in \chi(M)$ in the direction $V\in T_{p}M$. Since,
in general, the Chern connection coefficients $\Gamma^{i}_{jk}$ in
natural coordinates have a directional dependence, we must say
explicitly that $D_{V}U$ is defined with a fixed reference vector.
In particular, let $\sigma :[0,r]\longrightarrow M$ be a smooth
curve with velocity field $T=T(t)=\dot{\sigma}(t)$. Suppose that $U$
and $W$ are vector fields defined along $\sigma$. We define $D_{T}U$
with \emph{reference vector} $W$ as
$$D_{T}U=\left[\frac{dU^{i}}{dt}+U^{j}T^{k}(\Gamma^{i}_{jk})_{(\sigma ,W)}\right]\frac{\partial}{\partial x^{i}}\mid_{\sigma(t)}.$$
A curve $\sigma:[0,r]\longrightarrow M$, with velocity
$T=\dot{\sigma}$ is a Finslerian \emph{geodesic} if\\

$D_{T}\left[\frac{T}{F(T)}\right]=0$ ,\hskip.5cm with reference vector $T$.\\

We assume that all our geodesics $\sigma (t)$ have been
parameterized to have constant Finslerian speed. That is, the length
$F(T)$ is constant. These geodesics are characterized by the
equation \\

$D_{T}T=0$ , \hskip.5cm  with reference vector $T$.\\ \\Since
$T=\frac{d\sigma^{i}}{dt}\frac{\partial}{\partial x^{i}}$, this
equation says that
$$\frac{d^{2}\sigma^{i}}{dt^{2}}+\frac{d\sigma^{j}}{dt}\frac{d\sigma^{k}}{dt}(\Gamma^{i}_{jk})_{(\sigma ,T)}=0.$$

If $U,V$ and $W$ are vector fields along a curve $\sigma$, which has
velocity $T=\dot{\sigma}$, we have the derivative rule
$$\frac{d}{dt}g_{_W}(U,V)=g_{_W}(D_{T}U,V)+g_{_W}(U,D_{T}V)$$
whenever $D_{T}U$ and $D_{T}V$ are with reference vector $W$ and one
of the following conditions holds:
\begin{description}
    \item[i)] U or V is proportional to W, or
    \item[ii)] W=T and $\sigma$ is a geodesic.
\end{description}
\subsection{Homogeneous  Finsler Spaces}
\

Let $(M,F)$ be a Finsler space, where $F$ is positively homogeneous.
As in the Riemannian case, we have two kinds of definition of
isometry on $(M,F)$, in terms of Finsler function in the tangent
space and the induced non-reversible distance function on the base
manifold $M$. The equivalence of these two definitions in the
Finsler case is a result of S. Deng and Z. Hou [3]. They also prove
that the group of isometries of a Finsler space is a Lie
transformation group of the underlying manifold which can be used to
study homogeneous Finsler spaces.
\begin{Def}
A Finsler space $(M,F)$ is called  homogeneous Finsler space if the
group of isometries of $(M,F)$ , $I(M,F)$ , acts transitively on
$M$.
\end{Def}
A Finsler manifold $(M,F)$ is said to be forward geodesically
complete if every geodesic $\gamma(t)$, $a\leq t<b$, parameterized
to have constant Finslerian speed, can be extended to a geodesic
defined on $a\leq t <\infty$.
\begin{thm}[12]
Every  homogeneous Finsler space is forward complete.
\end{thm}
\begin{thm}[4]
Let $G$ be a Lie group, $H$ be a closed subgroup of $G$. Suppose
there exists an invariant Finsler metric on $\frac{G}{H}$. Then
there exists an invariant Riemannian metric on $\frac{G}{H}$.
\end{thm}
Let $M=\frac{G}{H}$ be a homogeneous space, where H is the isotropy
subgroup at a point $o\in M$. If the linear isotropy representation
$\lambda:H\longrightarrow GL(M_{o})$, $h\longrightarrow h_{\ast o}$,
is faithful, that is, injective, then $G$ acts effectively on $M$.
Let $(M,F)$ be a connected homogeneous Finsler manifold. If $G$ is
any connected transitive group of isometries of $M$ and $H$ is the
isotropy subgroup at a point, then $M$ is naturally identified with
the homogeneous manifold $\frac{G}{H}$. The Finsler metric $F$ on
$M$ can be considered as a G-invariant Finsler metric on
$\frac{G}{H}$. By Theorem 2.3, there exists a G-invariant Riemannian
metric on $\frac{G}{H}$. So the linear isotropy representation is
faithful and $G$ acts effectively on $\frac{G}{H}$.
 \\A homogeneous space $\frac{G}{H}$ is called reductive if there
exists a vector space decomposition
$\textbf{\underline{g}}=\textbf{m}\oplus \textbf{\underline{h}}$
such that $Ad(H)\textbf{m}\subset \textbf{m}$. In this case
$\textbf{m}\oplus \textbf{\underline{h}}$ is called a reductive
decomposition of $\textbf{\underline{g}}$. It is well-known
([7],[9]) that each Riemannian homogeneous space is reductive. We
now have the following.
\begin{rem}
Any homogeneous Finsler manifold $M=\frac{G}{H}$ is a reductive
homogeneous space.
\end{rem}

\section{Homogeneous geodesics in homogeneous Finsler
spaces} \

Let $(M=\frac{G}{H},F)$ be a homogeneous Finsler space with a fixed
origin $p$. Let $\underline{\textbf{g}}$ and
$\underline{\textbf{h}}$ be the Lie algebra of $G$ and $H$
respectively and let $$\underline{\textbf{g}}=\textbf{m}\oplus
\underline{\textbf{h}}$$ be a reductive decomposition of the Lie
algebra $\underline{\textbf{g}}$. From the Remark 2.4 such a
decomposition always exists.\\ For each $X\in
\underline{\textbf{g}}$ we obtain the corresponding fundamental
vector field $X^{\ast}$ on $M$ by means of
$$X^{\ast}_{q}=\frac{d}{dt}\mid_{t=0}\left(exp(tX)q\right)\hskip1cm \forall q\in M.$$
The canonical projection $\pi : G\longrightarrow \frac{G}{H}$
induces an isomorphism between the subspace $\textbf{m}$ and the
tangent space $T_{p}M$. Identifying  $\underline{\textbf{g}}$  with
$T_{e}G$  we get  $d\pi (X)=X^{\ast}_{p}$  for each $X\in
\underline{\textbf{g}}$ and hence $d\pi(X_{m})=X^{\ast}_{p}$.\\
Using this natural identification and scalar product
$g_{_{X^{\ast}_{p}}}$ on $T_{p}M$ we obtain a scalar product
$g_{_{X_{m}}}$ on $\textbf{m}$.

\

A vector $X\in \underline{\textbf{g}}-\{o\}$ will be called a
\emph{geodesic vector} if the curve $\gamma(t)=exp(tX)(p)$ is a
constant speed geodesic of $(M,F)$.\\ Let $(M=\frac{G}{H},g)$ be a
Riemannian homogeneous space, and
$\underline{\textbf{g}}=\textbf{m}\oplus \underline{\textbf{h}}$ be
a reductive decomposition. O. Kowalski and L. Vanhecke [10] proved
that $X\in \underline{\textbf{g}}$ is a geodesic vector if and only
if $$g([X,Y]_{m},X_{m})=0 \hskip.5cm \forall Y\in \textbf{m}.$$ In
the Finslerian case we get the following theorem. We use some ideas
from [10] in our proof.

\begin{thm}
A vector $X\in \textbf{\underline{\textsc{g}}}-\{0\}$ is geodesic
vector if and only if
$$g_{_{X_{m}}}(X_{m},[X,Z]_{m})=0\hskip.5cm \forall   Z\in
\textbf{\underline{\textsc{g}}}.$$
\end{thm}
 Proof: Let $(M,F)$ be a Finsler space. For any vector fields $T, V,
W$ on $M$, we have [2]
\begin{equation}
Tg_{_W}(V,W)=g_{_W}(D_{T}V,W)+g_{_W}(V,D_{T}W) \qquad \mbox{with
reference W}
\end{equation}
Similarly,
\begin{equation}
Vg_{_W}(T,W)=g_{_W}(D_{V}T,W)+g_{_W}(T,D_{_V}W),
\end{equation}
\begin{equation}
Wg_{_W}(V,W)=g_{_W}(D_{W}V,W)+g_{_W}(V,D_{_W}W).
\end{equation}
All covariant derivatives have $W$ as reference vector.\\
Subtracting (2) from the summation of (1) and (3) we get
\begin{eqnarray*}
  g_{_W}(V,D_{W+T}W)+g_{_W}(W-T,D_{_V}W) &=& Tg_{_W}(V,W)-Vg_{_W}(T,W)+Wg_{_W}(V,W) \\
  && -g_{_W}([T,V],W)-g_{_W}([W,V],W),
\end{eqnarray*}
where we have used the symmetry of the connection, i.e.,
$D_{V}W-D_{W}V=[V,W]$. Set $T=W-V$ in the above equation, we obtain
\begin{equation}
  2g_{_W}(V,D_{W}W)=2Wg_{_W}(V,W)-Vg_{_W}(W,W)-2g_{_W}([W,V],W).\\
\end{equation}
Let $X,Z \in \textbf{\underline{g}}$ be given and denote by $X^{*}$
and $Z^{*}$ the corresponding fundamental vector fields on $M$. From
the above equation we get
\begin{equation}
2g_{_X*}(D_{X^{*}}X^{*},Z^{*})=2X^{*}g_{_X*}(X^{*},Z^{*})-Z^{*}g_{_X*}(X^{*},X^{*})+2g_{_X*}([Z^{*},X^{*}],X^{*}).
\end{equation}
Recall also the formulas
\begin{equation}
Ad(exp (tX))Y= Y+t[X,Y]+O(t^{2}), \hskip1cm X,Y \in
\textbf{\underline{\textsc{g}}},
\end{equation}
\begin{equation}
k.exp(tX).k^{-1}=exp(tAd(k)X), \hskip1cm k\in G, X\in
\textbf{\underline{\textsc{g}}}.
\end{equation}
Denote briefly $g_{t}=exp(tX)$ , $h_{s}=exp(sZ)$. Using (7) and (6),
we get first, for any $x\in M$,
\begin{eqnarray*}
  Z^{*}_{g_{t}(x)}=\frac{d}{ds}\mid_{0}h_{s}g_{t}(x) &=& (dg_{t})\frac{d}{ds}\mid_{0}g_{t}^{-1}h_{s}g_{t}(x) \\
   &=& (dg_{t})\frac{d}{ds}\mid_{0}exp(s Ad(g_{t}^{-1})Z)(x) \\
   &=& (dg_{t})[Ad(g_{t}^{-1})Z]^{*}_{x}\\
   &=& (dg_{t})[Ad(exp(-tX))Z]^{*}_{x}\\
   &=& (dg_{t})[Z-t[X,Z]+O(t^{2})]_{x}^{*}.
\end{eqnarray*}
Similarly, we get
$$X^{*}_{h_{s}(x)}=(dh_{s})[X-s[Z,X]+O(s^{2})]_{x}^{*}.$$
We shall also use the obvious relations
$$X^{*}_{g_{t}(x)}=(dg_{t})X^{*}_{x} \hskip.2cm ,\hskip.2cm Z^{*}_{h_{s}(x)}=(dh_{s})Z^{*}_{x}.$$
Since $g_{t}$ is an isometry, $dg_{t}$ is a linear isometry between
the spaces $T_{p}M$ and $T_{g_{t}(p)}M$, $\forall p \in M$.
Therefore for any vector fields $V, W$ on $M$ we have
\begin{equation}
g_{dg_{t}(X^{*})}(dg_{t}(V),dg_{t}(W))=g_{_X*}(V,W).
\end{equation}
Now, we calculate
\begin{eqnarray*}
  X^{*}_{x}g_{_X*}(X^{*},Z^{*}) &=& \frac{d}{dt}\mid_{0}g_{_X*}(X^{*}_{g_{_t}(x)},Z^{*}_{g_{_t}(x)})  \\
   &=& \frac{d}{dt}\mid_{0}g_{_X*}(dg_{_t}(X^{*}_{x}),(dg_{_t})[Z-t[X,Z]+O(t^{2})]^{*}_{x}) \\
   &=& \frac{d}{dt}\mid_{0}g_{_X*}(X^{*}_{x},Z^{*}_{x}+t[X^{*},Z^{*}]_{x}+O(t^{2})) \\
   &= &g_{_X*}(X^{*},[X^{*},Z^{*}])(x).
\end{eqnarray*}
Further,
\begin{eqnarray*}
  Z^{*}_{x}g_{_X*}(X^{*},X^{*}) &=& \frac{d}{ds}\mid_{0}g_{_X*}(X^{*}_{x}+s[Z^{*},X^{*}]_{x}+O(s^{2}),X^{*}_{x}+s[Z^{*},X^{*}]_{x}+O(s^{2})) \\
   &=& 2g_{_X*}(X^{*},[Z^{*},X^{*}])(x).
\end{eqnarray*}
Substituting into (5), we get on $M$
\begin{equation}
g_{_X*}(D_{X^{*}}X^{*},Z^{*})=g_{_X*}(X^{*},[X^{*},Z^{*}])=-g_{_X*}(X^{*},[X,Z]^{*}).
\end{equation}
Now, suppose first that $X$ is a geodesic vector i.e.
$\gamma(t)=exp(tX)(p)$ is a geodesic of $F$ with Finslerian constant
speed. Then $D_{X^{*}_{g_{_t}(p)}} X^{*}_{g_{_t}(p)}=0$, so in
particular, $$g_{_X*}(X^{*}_{p},[X,Z]_{p}^{*})=0.$$ Using the
natural identification of $\textbf{m}$ and $T_{p}M$ we obtain
$$g_{_{X_{m}}}(X_{m},[X,Z]_{m})=0.$$

Let $Z$ be an arbitrary vector field on $M$ by (5) and (8) we have
\begin{eqnarray*}
   2g_{dg_{_t}(X^{*})}(dg_{_t}(Z),D_{dg_{_t}X^{*}}dg_{_t}X^{*}) &=&2(dg_{t}X^{\ast})g_{dg_{_t}(X^{*})}(dg_{_t}Z,dg_{_t}X^{*}) \\
   &&  -dg_{_t}(Z)g_{dg_{_t}(X^{*})}(dg_{_t}X^{*},dg_{_t}X^{*})\\
   &&  +2g_{dg_{_t}(X^{*})}([dg_{_t}Z,dg_{_t}X^{*}],dg_{_t}X^{*})\\
   &=& 2X^{\ast}g_{_X*}(Z,X^{*})-Zg_{_X*}(X^{*},X^{*})+2g_{_X*}([Z,X^{*}],X^{*}) \\
   &=& 2g_{_X*}(Z,D_{X^{*}}X^{*}).
\end{eqnarray*}
Consequently
$$g_{dg_{_t}(X^{*})}(dg_{_t}(Z),D_{dg_{_t}X^{*}}dg_{_t}X^{*})=g_{dg_{_t}(X^{*})}(dg_{_t}(Z),dg_{_t}(D_{X^{*}}X^{*})).$$
Since $Z$ is arbitrary and $g_{dg_{_t}(X^{*})}(.,.)$ is an inner
product, we have
\begin{equation}
D_{dg_{_t}X^{*}}dg_{_t}X^{*}=dg_{_t}(D_{X^{*}}X^{*}).
\end{equation}
On the other hand, suppose that
$g_{_{X^{*}_{p}}}(X^{*}_{p},[X,Z]_{p}^{*})=0$. Then
\begin{eqnarray*}
  g_{_{X_{g_{t}(p)}^{*}}}(X^{*},[X,Z]^{*})(exp(tX)(p)) &=& g_{_{X_{g_{t}(p)}^{*}}}(dg_{_t}X^{*}_{p},dg_{_t}[X,Z]_{p}^{*}) \\
   &=&g_{_{X^{*}_{p}}}(X^{*}_{p},[X,Z]^{*}_{p})=0
\end{eqnarray*}
for any $Z\in \textbf{\underline{g}}.$ Then (9) yields that
\begin{eqnarray*}
   g_{_{X^{*}}}(D_{X^{*}}X^{*},Z)(exp(tX)(p))&=& g_{_{X^{*}}}(dg_{_t}(D_{X^{*}}X^{*})_{p},dg_{_t}Z^{*}_{p}) \\
   &=& g_{_{X^{*}}}((D_{dg_{_t}X^{*}}dg_{_t}X^{*})_{p},dg_{_t}Z^{*}_{p}) \\
   &=& g_{_{X^{*}}}(D_{X^{*}_{g_{_t}(p)}}
   X^{*}_{g_{_t}(p)},dg_{_t}Z^{*}_{p})=0.
\end{eqnarray*}
Then this yields that $exp(tX)(p)$ is a geodesic with constant
speed.\\\qed
\begin{cor}
A vector $X \in \textbf{\underline{\textsc{g}}}-\{0\}$ is geodesic
vector  if and only if
\begin{equation}
 g_{_{X_{m}}}(X_{m},[X,Y]_{m})=0 \qquad \mbox{for all  $Y\in \textbf{m}$.}
\end{equation}
\end{cor}
Proof: Since $F$ is G-invariant, we have
$$F(Ad(h)W)=F(W)\hskip.5cm \forall h\in H , W\in \textbf{m}.$$
Therefore, $\forall y\neq 0 , u , v \in \textbf{m}$, $x\in
\underline{\textbf{h}},t , r, s\in R$, we have
$$F^{2}(Ad(exp(tx)))(y+ru+sv)=F^{2}(y+ru+sv).$$
By definition,
$$g_{y}(u,v)=\frac{1}{2}\frac{\partial^{2}}{\partial r \partial s}F^{2}(y+ru+sv)\mid_{r=s=0}.$$
Thus $$g_{y}(u,v)=\frac{1}{2}\frac{\partial^{2}}{\partial r \partial
s}F^{2}(Ad(exp(tx))(y+ru+sv))\mid_{r=s=0}.$$ Now for $w\in
\textbf{m}$, from (6) we have
$$Ad(exp(tx))w =w+t[x,w]+O(t^{2}).$$ Therefore
$$g_{y}(u,v)=\frac{1}{2}\frac{\partial^{2}}{\partial r
\partial s}F^{2}(y+ru+sv+t[x,y+ru+sv]+O(t^{2}))\mid_{r=s=0}.$$
Taking derivative with respect to $t$ at $t=0$, we get
\begin{equation}
0=g_{y}([x,u],v)+g_{y}(u,[x,v])+2C_{y}([x,y],u,v),
\end{equation}
where $C_{y}$ is the Cartan tensor of $F$ at $y$. It follows from
the homogeneity of $F$ that $C_{y}(y,v,w)=0$. So we have
$$g_{y}([x,y],y)=0.$$
For any $Z\in \textbf{\underline{g}}$, where $Z=Y+A$ with $Y\in
\textbf{m}$, $A\in \underline{\textbf{h}}$, we obtain
$$g_{_{X_{m}}}([X,Z]_{m},X_{m})=g_{_{X_{m}}}([X,Y]_{m},X_{m})+g_{_{X_{m}}}([X,A]_{m},X_{m}).$$
Here, the second term is equal to $g_{_{X_{m}}}([X_{m},A],X_{m})=0.$
Hence (11) implies that $X$ is a geodesic vector.\\\qed
\begin{cor}
If $X \in \textbf{\underline{\textsc{g}}}-\{0\}$ is a geodesic
vector then $Ad(h)X$ and $\lambda X$ are geodesic vector for all
$h\in H$, $\lambda \in R.$
\end{cor}
Proof: Evident from the fact that
$$g_{y}(u,v)=g_{_{Ad(h)y}}(Ad(h)u,Ad(h)v)\hskip.5cm \forall h \in H.$$\\\qed

 In the following theorem, we consider bi-invariant Finsler
metrics on Lie groups. We show that the geodesics of $G$ starting at
the identity element are the one-parameter subgroup of $G$. Let $G$
be a connected Lie group. S. Deng and Z. Hou [4] prove that there
exists a bi-invariant Finsler metric on $G$ if and only if there
exists a Minkowski norm $F$ on $\textbf{\underline{g}}$ such that
$$g_{y}([x,u],v)+g_{y}(u,[x,v])+2C_{y}([x,y],u,v)=0$$
$\forall y\in \underline{\textbf{\textsc{g}}}-\{0\},x,u,v\in
\textbf{\underline{g}}$. So we have the following:
\begin{thm}
Let $G$ be a connected Lie group furnished with a bi-invariant
Finsler metric $F$. Then each vector of
$\textbf{\underline{\textsc{g}}}$ is geodesic vector.
\end{thm}

Here we study the existence of homogeneous geodesics in homogeneous
Finsler spaces. The problem of the existence of homogeneous
geodesics in homogeneous Finsler manifolds seems to be an
interesting one. About the existence of homogeneous geodesics in a
general homogeneous Riemannian manifold, we have, at first, a result
due to V. V. Kajzar who proved that a Lie group endowed with a
left-invariant metric admits at least one homogeneous geodesic [6].
More recently O. Kowalski and J. Szenthe extended this result to all
homogeneous Riemannian manifolds [9]. Homogeneous geodesics of
left-invariant Lagrangian on Lie groups were studied by J. Szenthe
[20]. The following result is due to J. Szenthe [20].
\begin{thm}
Let $G$ be a compact connected Lie group and $L: TG\longrightarrow
R$ a left-invariant Lagrangian which is a first integral of its
Lagrangian field. Then $L$ has at least one homogeneous geodesic.
If, in particular, $G$ is also semi-simple and of rank $\geq 2$ then
$L$ has infinitely many homogeneous geodesics.
\end{thm}

Let $(M,F)$ be a Finsler space. For every smooth parameterized curve
$\gamma : [0,1]\longrightarrow M$, the length of $\gamma$ is given
by
\begin{equation}
L(\gamma)=\int_{0}^{1}F(\gamma(t),\dot{\gamma}(t))dt.
\end{equation}

A geodesic of the Finsler space $(M,F)$ is an extermal curve of
(13). This is in fact a solution of the Euler-Lagrange equations
\begin{equation}
\frac{d}{dt}(\frac{\partial F}{\partial \dot{x}^{i}})-\frac{\partial
F}{\partial x^{i}}=0, \hskip.5cm \dot{x}^{i}=\frac{dx^{i}}{dt}
\end{equation}
where $(x^{i}(t))$ is a local coordinate expression of $\gamma$.
This system is equivalent to
\begin{equation}
\frac{d^{2}x^{i}}{dt^{2}}+2G^{i}(x,\frac{dx}{dt})=0
\end{equation}
where $$G^{j}(y)=\frac{1}{4}g^{jl}(y)\left[2\frac{\partial
g_{sl}}{\partial x^{k}}(y)-\frac{\partial g_{sk}}{\partial
x^{l}}(y)\right] y^{s}y^{k}.$$ Let $$G=y^{i}\frac{\partial}{\partial
x^{i}}-2G^{i}(x,y)\frac{\partial}{\partial y^{i}}$$ $G$ is a vector
field on $TM-\{0\}$. It is easy to see that $x(t)$ is a solution of
(15) if and only if its lift $\dot{x}(t)=(x(t), \frac{dx}{dt}(t))$
is an integral curve of $G$ in $TM-\{0\}$. $G$ is called the
geodesic spray. The following lemma show that any Finsler metric $F$
is a first integral of its geodesic spray.
\begin{lem}
[17] For any Finsler metric $F$ on a manifold, $G(F)=0$
\end{lem}
So from Theorem 3.5 the following result follows:
\begin{thm}
Let $G$ be a compact connected Lie group and $F$ a left-invariant
Finsler metric. Then $F$ has at least one homogeneous geodesic. If,
in particular, $G$ is also semi-simple and of rank $\geq 2$ then $F$
has infinitely many homogeneous geodesics.
\end{thm}
\subsection{Naturally Reductive Homogeneous Finsler Space}
   \ \ \ The scheme is to treat the geometry of coset manifolds
   $\frac{G}{H}$ as a generalization of the geometry of Lie group
   $G$ ( Since $\frac{G}{H}$ reduces to $G$ when H=\{e\} ). From
   this viewpoint, the isomorphism $\textbf{m}\simeq T_{o}(\frac{G}{H})$
   generalizes the canonical isomorphism $\textbf{\underline{g}}\simeq T_{e}G$, and a
   G-invariant Riemannian metric on $\frac{G}{H}$ generalizes a
   left-invariant metric on $G$. The notion of bi-invariant
   Riemannian metric on $G$ generalizes as follows.
\begin{Def}
A Riemannian homogeneous space $(\frac{G}{H},g)$ is said to be
naturally reductive if there exists a reductive decomposition
$\textbf{\underline{\textsc{g}}}=\textbf{m}+\textbf{\underline{h}}$
of $\textbf{\underline{\textsc{g}}}$ satisfying the condition
\begin{equation}
<[X,Y]_{m},Z>+<Y,[X,Z]_{m}>=0
\end{equation}
for all $X,Y,Z \in \textbf{m}.$
\end{Def}
where $<,>$ denotes the inner product on $\textbf{m}$ induced by the
metric $g$.
\\
In fact, when $H=\{e\}$, hence $\textbf{m}=\textbf{\underline{g}}$,
the above condition is just the condition \begin{equation}
 <[X,Y],Z>+<Y,[X,Z]>=0,
\end{equation}
for a bi-invariant Riemannian metric on $G$. In [4] authors
introduced the notion of a Minkowski Lie algebra:
\begin{Def}
Let $\textbf{\underline{\textsc{g}}}$ be a real Lie algebra, $F$ be
a Minkowski norm on $\underline{\textbf{\textsc{g}}}$. Then
$\{\underline{\textbf{\textsc{g}}},F\}$ is called a Minkowski Lie
algebra if the following condition is satisfied
\begin{equation}
g_{y}([x,u],v)+g_{y}(u,[x,v])+2C_{y}([x,y],u,v)=0
\end{equation}
where $y\in \textbf{\underline{\textsc{g}}}-\{0\},x,u,v \in
\textbf{\underline{\textsc{g}}}$.
\end{Def}
They showed that
\begin{thm}
Let $G$ be a connected Lie group. Then there exists a bi-invariant
Finsler metric on $G$ if and only if there exists a Minkowski norm
$F$ on $\textbf{\underline{\textsc{g}}}$ such that
$\{\textbf{\underline{\textsc{g}}},F\}$ is a Minkowski Lie algebra.
\end{thm}
It is easy to see that the notion of Minkowski Lie algebra is the
natural generalization of (17). Now we define the notion of
naturally reductive homogeneous Finsler space.
\begin{Def}
A homogeneous manifold $\frac{G}{H}$ with an invariant Finsler
metric $F$ is called naturally reductive if there exists an
Ad(H)-invariant decomposition
$\textbf{\underline{\textsc{g}}}=\textbf{\underline{h}}+\textbf{m}$
such that
$$g_{y}([x,u]_{m},v)+g_{y}(u,[x,v]_{m})+2C_{y}([x,y]_{m},u,v)=0$$
where $y\neq 0,x,u,v \in \textbf{m}$.
\end{Def}
Evidently this definition is the natural generalization of (16). On
the other hand, when $H=\{e\}$, hence
$\textbf{m}=\textbf{\underline{g}}$, this formula is just the (18).
The following theorem is a consequence of Theorem 3.1 and the
Definition 3.11.
\begin{thm}
Let $(\frac{G}{H},F)$ be a naturally reductive homogeneous Finsler
space. Then each geodesic of $(\frac{G}{H},F)$ is an orbit of a
one-parameter group of isometries $\{ exp(tX) \}$, $X\in
\textbf{\underline{\textsc{g}}}$.
\end{thm}
\section{Homogeneous geodesics of Randers spaces}
\

In this section, we consider homogeneous geodesics in a  homogeneous
Randers space. Randers metrics were introduced by Randers in 1941
[16] in the context of general relativity. They are Finsler spaces
built from
\begin{description}
    \item[i)] a Riemannian metric $\widetilde{a}=\widetilde{a}_{ij}dx^{i}\otimes
    dx^{j}$, and
    \item[ii)] a 1-form $\widetilde{b}=\widetilde{b}_{i}dx^{i},$
\end{description}
both living globally on the smooth n-dimensional manifold $M$. The
Finsler function of a Randers metric has the simple form $F=\alpha+
\beta$, where
$$\alpha(x,y)=\sqrt{\widetilde{a}_{ij}(x)y^{i}y^{j}}\hskip.3cm,\hskip.3cm\beta(x,y)=\widetilde{b}_{i}(x)y^{i}.$$
 Generic  Randers metric are only positively homogeneous. No Randers
 metric can satisfy absolut homogeneity $F(x,cy)=|c|F(x,y)$ unless
 $\widetilde{b}=0$, in which case it is Riemannian. Also, in order
 for $F$ to be positive and strongly convex on $TM \backslash\{0\}$, it is
 necessary and sufficient to have $$\|\widetilde{b}\|=\sqrt{\widetilde{b}_{i}\widetilde{b}^{i}}<1,\hskip1cm where \hskip.5cm \widetilde{b}^{i}=\widetilde{a}^{ij}\widetilde{b}_{j}.$$
See [2]. Strong convexity means that the fundamental tensor $g_{ij}$
is positive definite. The Riemannian metric
$\widetilde{a}=\widetilde{a}_{ij}dx^{i}\otimes dx^{j}$ induces the
musical bijections between 1-forms and vector fields on $M$, namely
$\flat:T_{x}M\longrightarrow T_{x}^{\ast}M$ given by
$y\longrightarrow \widetilde{a}_{x}(y,\circ)$ and its inverse
$\sharp:T_{x}^{\ast}M\longrightarrow T_{x}M $. In the local
coordinates we have
$$(y^{\flat})_{i}=\widetilde{a}_{ij}y^{j}\hskip.5cm (\theta^{\sharp})^{i}=\widetilde{a}^{ij}\theta_{j}\hskip.5cm y\in T_{x}M\hskip.3cm \theta\in T_{x}^{\ast}M.$$
Now the corresponding vector field to the 1-form $\widetilde{b}$
will be denoted by $\widetilde{b}^{\sharp}$, obviously we have
$\|\widetilde{b}\|=\|\widetilde{b}^{\sharp}\|$ and
$$\beta(x,y)=(\widetilde{b}^{\sharp})^{\flat}(y)=\widetilde{a}_{x}(\widetilde{b}^{\sharp},y).$$
Thus a Randers metric $F$ with Riemannian metric
$\widetilde{a}=\widetilde{a}_{ij}dx^{i}\otimes dx^{j}$ and 1-form
$\widetilde{b}$ can be showed by
$$F(x,y)=\sqrt{\widetilde{a}_{x}(y,y)}+\widetilde{a}_{x}(\widetilde{b}^{\sharp},y)\hskip.5cm x\in M , y\in T_{x}M,$$
where
$\widetilde{a}_{x}(\widetilde{b}^{\sharp},\widetilde{b}^{\sharp})<1,
\hskip.5cm \forall x\in M$.
\begin{thm}
Let $(M,F)$ be a homogeneous Randers space with $F$ defined by the
Riemannian metric $\widetilde{a}$ and the vector field $X$. Then $X$
is a geodesic vector of $(M,\widetilde{a})$ if and only if $X$ is a
geodesic vector of $(M,F)$.
\end{thm}
Proof: Let
$F(p,y)=\sqrt{\widetilde{a}_{p}(y,y)}+\widetilde{a}_{p}(X,y)$.\\ Now
for $s,t \in R$
\begin{eqnarray*}
  F^{2}(y+su+tv) &=& \widetilde{a}(y+su+tv,y+su+tv)+\widetilde{a}^{2}(X,y+su+tv) \\
   & & +2\sqrt{\widetilde{a}(y+su+tv,y+su+tv)} \widetilde{a}(X,y+su+tv)
\end{eqnarray*}
By definition
$$g_{y}(u,v)=\frac{1}{2}\frac{\partial^{2}}{\partial r \partial
s}F^{2}(y+ru+sv)\mid_{r=s=0}.$$ So by a direct computation we get

\begin{eqnarray*}
  g_{y}(u,v) &=& \widetilde{a}(u,v)+
  \widetilde{a}(X,u)\widetilde{a}(X,v)\\\\
  && +\frac{\widetilde{a}(u,v)\widetilde{a}(X,y)}{\sqrt{\widetilde{a}(y,y)}}-\frac{\widetilde{a}(v,y)\widetilde{a}(u,y)\widetilde{a}(X,y)}{\widetilde{a}(y,y)\sqrt{\widetilde{a}(y,y)}} \mbox{\qquad
  (19)}\\\\
   &&+\frac{\widetilde{a}(X,v)\widetilde{a}(u,y)}{\sqrt{\widetilde{a}(y,y)}}+\frac{\widetilde{a}(X,u)\widetilde{a}(v,y)}{\sqrt{\widetilde{a}(y,y)}}.\\
\end{eqnarray*}
So for all $Z\in \textbf{m}$ we have
\begin{eqnarray*}
  g_{_X}(X,[X,Z]_{m}) &=& \widetilde{a}(X,[X,Z]_{m})+\widetilde{a}(X,X)\widetilde{a}(X,[X,Z]_{m}) \\
   && +2\sqrt{\widetilde{a}(X,X)}\widetilde{a}(X,[X,Z]_{m}) \\
   &=&
   \widetilde{a}(X,[X,Z]_{m})(1+\sqrt{\widetilde{a}(X,X)}+F(X)).
\end{eqnarray*}
Thus $g_{_X}(X,[X,Z]_{m})=0$ if and only if
  $\widetilde{a}(X,[X,Z]_{m})=0$.\\
\qed
\begin{thm}
Let $(M,F)$ be a homogeneous Randers space with $F$ defined by the
Riemannian metric $\widetilde{a}$ and the vector field $X$. Let
$y\in \textbf{\underline{\textsc{g}}}$ be a vector which
$\widetilde{a}(X,[y,z]_{m})=0$ for all $z \in \textbf{m}$. Then $y$
is a geodesic vector of $(M,F)$ if and only if $y$ is a geodesic
vector of $(M,\widetilde{a})$.
\end{thm}
Proof: According to the above formula for $g_{y}(u,v)$, we have
\begin{eqnarray*}
  g_{y_{_m}}(y_{m},[y,z]_{m}) &=& \widetilde{a}(y_{m},[y,z]_{m})+\widetilde{a}(X,y_{m})\widetilde{a}(X,[y,z]_{m})
  \\ \\
   &&
   +\frac{\widetilde{a}(y_{m},[y,z]_{m})\widetilde{a}(X,y_{m})}{\sqrt{\widetilde{a}(y_{m},y_{m})}}+\widetilde{a}(X,[y,z]_{m})\sqrt{\widetilde{a}(y_{m},y_{m})}\\\\
   &=&\widetilde{a}(y_{m},[y,z]_{m})\left(1+\frac{\widetilde{a}(X,y_{m})}{\sqrt{\widetilde{a}(y_{m},y_{m})}}\right)+\widetilde{a}(X,[y,z]_{m})\left(\widetilde{a}(X,y_{m})+\sqrt{\widetilde{a}(y_{m},y_{m})}\right).
\end{eqnarray*}
So we have
\begin{eqnarray*}
  g_{y_{_m}}(y_{m},[y,z]_{m}) &=& \widetilde{a}(y_{m},[y,z]_{m})\left(\frac{F(y_{m})}{\sqrt{\widetilde{a}(y_{m},y_{m})}}\right) \\
   && +\widetilde{a}(X,[y,z]_{m})F(y_{m})
\end{eqnarray*}
This conclude the proof. \\\qed \\ Let $(M,F)$ be a Finsler space.
Then $(M,F)$ is called a Berwald space if the Chern connection
coefficients $\Gamma_{ij}^{k}(x,y)$ in natural coordinate systems
have no dependence on the vector $y$, or in other words, if the
Chern connection defined a linear connection directly on the
underlying manifold.
\begin{thm}
Let $(M,F)$ be a homogeneous Randers space with $F$ defined by the
Riemannian metric $\widetilde{a}=\widetilde{a}_{ij}dx^{i}\otimes
dx^{j}$ and the vector field $X$ which is of Berwald type. Then
$(M,F)$ is naturally reductive if and only if the underlying
Riemannian metric $(M,\widetilde{a})$ is naturally reductive.
\end{thm}
Proof: Let $(M,\widetilde{a})$ is naturally reductive. We show that
for all $0\neq y ,z,u,v \in \textbf{m}$
$$g_{y}([z,u]_{m},v)+g_{y}(u,[z,v]_{m})+2C_{y}([z,y]_{m},u,v)=0.$$
Since $F$ is of Berwald type, $(M,F)$ and $(M,\widetilde{a})$ have
the same connection. So according to the relation
\begin{eqnarray*}
  g_{y_{_m}}(y_{m},[y,z]_{m}) &=& \widetilde{a}(y_{m},[y,z]_{m})\left(\frac{F(y_{m})}{\sqrt{\widetilde{a}(y_{m},y_{m})}}\right) \\
   && +\widetilde{a}(X,[y,z]_{m})F(y_{m}),
\end{eqnarray*}
for all $0\neq y\in\textbf{m}$ we have
$$\widetilde{a}(X,[y,z]_{m})=0\hskip.4cm \forall z \in \textbf{m}.$$
From (19) we get
\begin{eqnarray*}
  g_{y}([z,u]_{m},v) &=& \widetilde{a}([z,u]_{m},v)\left(1+\frac{\widetilde{a}(X,y)}{\sqrt{\widetilde{a}(y,y)}}\right)\\
  && +\widetilde{a}([z,u]_{m},y)\left( \frac{\widetilde{a}(X,v)}{\sqrt{\widetilde{a}(y,y)}}-\frac{\widetilde{a}(v,y)\widetilde{a}(X,y)}{\widetilde{a}(y,y)\sqrt{\widetilde{a}(y,y)}}\right)
\end{eqnarray*}
\begin{eqnarray*}
  g_{y}(u,[z,v]_{m}) &=&\widetilde{a}(u,[z,v]_{m})\left(1+\frac{\widetilde{a}(X,y)}{\sqrt{\widetilde{a}(y,y)}}\right)  \\
   &&
   +\widetilde{a}([z,v]_{m},y)\left(\frac{\widetilde{a}(X,u)}{\sqrt{\widetilde{a}(y,y)}}-\frac{\widetilde{a}(u,y)\widetilde{a}(X,y)}{\widetilde{a}(y,y)\sqrt{\widetilde{a}(y,y)}}\right).
\end{eqnarray*}
By definition
$$C_{y}(z,u,v)=\frac{1}{2}\frac{d}{dt}[g_{y+tv}(z,u)]|_{t=0}.$$
So by a direct copmutation we get\\\\
\begin{eqnarray*}
  2C_{y}(z,u,v) &=&\frac{\widetilde{a}(z,u)\widetilde{a}(X,v)\widetilde{a}(y,y)-\widetilde{a}(y,v)\widetilde{a}(z,u)\widetilde{a}(X,y)}{\widetilde{a}(y,y)\sqrt{\widetilde{a}(y,y)}}  \\
   && +\frac{\widetilde{a}(X,u)\widetilde{a}(z,v)\widetilde{a}(y,y)-\widetilde{a}(y,v)\widetilde{a}(X,u)\widetilde{a}(z,y)}{\widetilde{a}(y,y)\sqrt{\widetilde{a}(y,y)}} \\
   &&  + \frac{\widetilde{a}(X,z)\widetilde{a}(u,v)\widetilde{a}(y,y)-\widetilde{a}(y,v)\widetilde{a}(X,z)\widetilde{a}(u,z)}{\widetilde{a}(y,y)\sqrt{\widetilde{a}(y,y)}}\\
   &&  -\frac{\widetilde{a}(x,y)\widetilde{a}(u,y)\widetilde{a}(z,v)+\widetilde{a}(x,y)\widetilde{a}(u,v)\widetilde{a}(z,y)+\widetilde{a}(x,v)\widetilde{a}(u,y)\widetilde{a}(z,y)}{\widetilde{a}(y,y)\sqrt{\widetilde{a}(y,y)}}\\
   &&
   +\frac{3\widetilde{a}(y,v)\widetilde{a}(X,y)\widetilde{a}(u,y)\widetilde{a}(z,y)}{\widetilde{a}(y,y)^{2}\sqrt{\widetilde{a}(y,y)}}.
\end{eqnarray*}
So we have
\begin{eqnarray*}
  C_{y}([z,y]_{m},u,v)
  &=&\widetilde{a}([z,y]_{m},u)\left(\frac{\widetilde{a}(X,v)}{\sqrt{\widetilde{a}(y,y)}}-\frac{\widetilde{a}(y,v)\widetilde{a}(X,y)}{\widetilde{a}(y,y)\sqrt{\widetilde{a}(y,y)}}\right)
  \\  \\
   &&+\widetilde{a}([z,y]_{m},v)\left(\frac{\widetilde{a}(X,u)}{\sqrt{\widetilde{a}(y,y)}}-\frac{\widetilde{a}(X,y)\widetilde{a}(u,y)}{\widetilde{a}(y,y)\sqrt{\widetilde{a}(y,y)}}\right).
\end{eqnarray*}
Therefore
\begin{eqnarray*}
  g_{y}([z,u]_{m},v)+g_{y}(u,[z,v]_{m})+2C_{y}([z,y]_{m},u,v)= &&  \\ \\
   (\widetilde{a}([z,u]_{m},v)+\widetilde{a}(u,[z,v]_{m}))\left(1+\frac{\widetilde{a}(X,y)}{\sqrt{\widetilde{a}(y,y)}}\right)&&  \\ \\
   +(\widetilde{a}([z,u]_{m},y)+\widetilde{a}([z,y]_{m},u))\left(\frac{\widetilde{a}(X,v)}{\widetilde{a}(y,y)}-\frac{\widetilde{a}(v,y)\widetilde{a}(X,y)}{\widetilde{a}(y,y)\sqrt{\widetilde{a}(y,y)}}\right)&&  \\ \\
   +(\widetilde{a}([z,v]_{m},y)+\widetilde{a}([z,y]_{m},v))\left(\frac{\widetilde{a}(X,u)}{\sqrt{\widetilde{a}(y,y)}}-\frac{\widetilde{a}(X,y)\widetilde{a}(u,y)}{\widetilde{a}(y,y)\sqrt{\widetilde{a}(y,y)}}\right).&&
\end{eqnarray*}
Thus
$$g_{y}([z,u]_{m},v)+g_{y}(u,[z,v]_{m})+2C_{y}([z,y]_{m},u,v)=0$$
for all $y\neq0,u,v,z \in \textbf{m}$.\

\

Conversely Let $(M,F)$ be naturally reductive i.e. for all
$y\neq0,u,v,z \in \textbf{m}$
$$g_{y}([z,u]_{m},v)+g_{y}(u,[z,v]_{m})+2C_{y}([z,y]_{m},u,v)=0.$$
So $$g_{y}([y,u]_{m},v)+g_{y}(u,[y,v]_{m})=0.$$ $(M,F)$ and
$(M,\widetilde{a})$ have the same geodesics and for all $0\neq y\in
\textbf{m}$, $y$ is a geodesic vector, so for all $y\in \textbf{m}$
we have $$\widetilde{a}(X,[y,z]_{m})=0 \hskip.4cm \forall z\in
\textbf{m}.$$ Therefore
$$g_{y}([y,u]_{m},v)=\widetilde{a}([y,u]_{m},v)\left(1+\frac{\widetilde{a}(X,y)}{\sqrt{\widetilde{a}(y,y)}}\right),$$
$$g_{y}(u,[y,v]_{m})=\widetilde{a}(u,[y,v]_{m})\left(1+\frac{\widetilde{a}(X,y)}{\sqrt{\widetilde{a}(y,y)}}\right).$$
So we have
$$g_{y}([y,u]_{m},v)+g_{y}(u,[y,v]_{m})=(\widetilde{a}([y,u]_{m},v)+\widetilde{a}(u,[y,v]_{m}))\left(1+\frac{\widetilde{a}(X,y)}{\sqrt{\widetilde{a}(y,y)}}\right)=0.$$
We easily see that
$\left(1+\frac{\widetilde{a}(X,y)}{\sqrt{\widetilde{a}(y,y)}}\right)\neq
0.$ Thus
$$\widetilde{a}([y,u]_{m},v)+\widetilde{a}(u,[y,v]_{m})=0$$\qed

\section{Some Curvature Properties}
\

The \textbf{S}-curvature is one of most important non-Riemannian
quantities in Finsler geometry which vanishes for Riemannian
metrics. The \textbf{S}-curvature has been introduced in [18]. In
this section, we discus the relationship between the homogeneous
geodesics and \textbf{S}-curvature.

Let $F$ be a Finsler metric on a manifold $M$. Let
$\{e_{i}\}_{i=1}^{n}$ be a basis for $T_{x}M$ and $\{\omega
_{i}\}_{i=1}^{n}$ the dual basis for $T_{x}^{^*}M$. Denote by
$d\mu_{x}= \sigma(x)\omega^{1}\wedge...\wedge\omega^{n}$ the
Busemann volume form at $x$, where
$$\sigma(x)=\frac{Vol(B^{n})}{Vol\{(y^{i})\in R^{n}, F(y^{i}e_{i})<1\}} ,$$
where $B^{n}$ denotes the unit ball in $R^{n}$ an $Vol$ denotes the
Euclidean measure on $R^{n}$. For each $y\in T_{x}M-\{0\}$, define
$$\tau(x,y)=Ln\left[\frac{\sqrt{detg_{y}(e_{i},e_{j})}}{\sigma(x)}\right].$$
The scalar function $\tau : TM\setminus
 \{0\} \longrightarrow R$ is called the \emph{distortion}. To
 measure the rate of changes of the distortion along geodesics, we
 define $$\textbf{S}(x,y)=\frac{d}{dt}\left[\tau (c(t),\dot{c}(t))\right]_{t=0}$$
 where $c(t)$ is the geodesic with $\dot{c}(0)=y$. \\ The scalar
 function $\textbf{S} : TM\setminus
 \{0\} \longrightarrow R$ is called the \textbf{S}-curvature ([17],[18]).
\begin{thm}
Let $X$ be a geodesic vector. Then $\textbf{S}(X_{m})=0$ for the
Busemann volume form.
\end{thm}
Proof: Let $\gamma(t)= exp(tX)(p)$ be the homogeneous geodesic
corresponding to $X$. Denote $g_{t}=exp(tX)$, obviously we have
$$X^{*}_{g_{_t}(x)}=(dg_{_t})(X^{*}_{x}).$$ Take an arbitrary basis
$\{e_{i}\}_{i=1}^{n}$ for $T_{x}M$. we obtain a frame along
$\gamma(t)$, $$e_{i}(t)=(dg_{_t})e_{i} \hskip.5cm i=1,...,n.$$ Let
$d\mu$ denote the Busemann volume form of $F$. Put $$d\mu
|_{c(t)}=\sigma(t)\omega^{1}(t)\wedge...\wedge\omega^{n}(t),$$ where
$\{\omega^{i}(t)\}$ be the basis for $T^{^\ast}_{\gamma(t)}M$ which
are dual to $\{e_{i}(t)\}_{i=1}^{n}.$
$$\sigma(t)=\frac{Vol(B^{n})}{Vol\{(y^{i})\in R^{n}, F(y^{i}e_{i}(t))<1\}}.$$
Since $g_{_t}$ is an isometry we have
$F(y^{i}e_{i}(t))=F(y^{i}e_{i})$ and by the relation (8) we know
that
$$g_{\dot{\gamma}(t)}(e_{i}(t),e_{j}(t))=g_{\dot{\gamma}(0)}(e_{i},e_{j}).$$
Thus
$$det\left[g_{\dot{\gamma}(t)}(e_{i}(t),e_{j}(t))]=det[g_{\dot{\gamma}(0)}(e_{i},e_{j})\right],$$
$$\{(y^{i})\in R^{n}, F(y^{i}e_{i}(t))<1\}= \{(y^{i})\in R^{n}, F(y^{i}e_{i})<1\}.$$
This implies
$$\tau(\gamma(t),\dot{\gamma}(t))=Ln\left[\frac{\sqrt{det[g_{\dot{\gamma}(t)}(e_{i}(t),e_{j}(t))]}}{\sigma(t)}\right]$$
is constant, so $\textbf{S}(\gamma(t),\dot{\gamma}(t))=0$. In
particular, at $t=0$, $\textbf{S}(X^{*}_{p})=0$.\\\qed
\\
Now similar to Riemannian spaces we define the notion of
\emph{geodesic orbit (g.o)} space for Finsler spaces.
\begin{Def}
A homogeneous Finsler space $(M,F)$ is said to be geodesic orbit
(g.o) space if every geodesic in $M$ is an orbit of a one-parameter
group of isometries i.e. there exists a transitive group $G$ of
isometries such that every geodesic in $M$ is of the form $exp(tX)p$
with $X\in \textbf{\underline{\textsc{g}}}$ , $p\in M$.
\end{Def}
\begin{cor}
Let $(M,F)$ be a g.o Finsler space. Then the \textbf{S}-curvature
\textbf{S}=0 for the Busemann volume form.
\end{cor}
Proof: Let $0\neq X \in T_{x}M $ be an arbitrary vector. Let
$\gamma$ be the geodesic with $\dot{\gamma}(0)=X$. Since $(M,F)$ is
g.o we can write $$\gamma(t)=exp(tX)\gamma(0).$$ According to
Theorem 5.1, $\textbf{S}(\gamma(t),\dot{\gamma}(t))=0$ and so
\textbf{S}(x,X)=0.\qed
\section*{Acknowledgements}
 I am grateful to the referee for the valuable suggestions and
 comments .\`{a}Çð?

\noindent


\begin{thebibliography}{99}




\bibitem[1]{1}V. I. Arnold, Sur la g\'{e}om\'{e}trie
diff\'{e}rentielle des groupes de Lie de dimension infinie et ses
applications \`{a} l'hydrodynamique des fluides parfaites, Ann.
Inst. Fourier(Grenoble), 16(1960), 319-361.



\bibitem[2]{2} D.Bao, S.S. Chern and Shen, An Introduction to
Riemann-Finsler geometry, Springer-Verlag,New-York.2000.


\bibitem[3]{3}S. Deng and Z. Hou, The group of isometries of a
Finsler space, Pacific.J.Math 207(1) (2002), 149-155.

\bibitem[4]{4}S. Deng and Z. Hou, Invariant Finsler metrics on
homogeneous manifolds, J.Phys.A: Math. Gen. 37(2004) 8245-8253.
\bibitem[5]{5}C. Gordon, Homogeneous Riemannian manifolds whose
geodesics are orbits, Prog. Nonlinear Differential Equations Appl.
20 (1996) 155-174.
\bibitem[6]{6}V.V. Kajzer, Conjugate points of left-invariant
metrics on Lie groups, Soviet Math. 34(1990), 32-44.
\bibitem[7]{7}S. Kobayashi and K. Nomizu, Foundations of
Differential Geometry II, Interscience Publishers, New York, 1969.
\bibitem[8]{8}O. Kowalski, S. Nik\u{c}evi\'{c} and Z.
Vl\'{a}\u{s}ek, Homogeneous geodesics in homogeneous Riemannian
manifolds-Examples, Geometry and Topology of Submanifolds
(Beijing/Berlin 1999)(2000),World Sci. Publishing Co., River Edge,
NJ, 104-112.
\bibitem[9]{9}O. Kowalski and J. Szenthe, On the existence of
homogeneous geodesics in homogeneous Riemannian manifolds, Geom.
Dedicata, 81(2000), 209-214. Erratum: Geom. Dedicata, 84 (2001),
331-332.
\bibitem[10]{10}O. Kowalski and L. Vanhecke, Riemannian manifolds with
homogeneous geodesics, Boll. Un. Mat. Ital., 5(1991), 189-246.
\bibitem[11]{11}A. Krist\'{a}ly and L. Kozma, Metric characterization of Berwald spaces of
non-positive flag curvature, J. Geom. Phys., 56(2006), 1257-1270.


\bibitem[12]{12}D. Latifi and A. Razavi, On homogeneous Finsler
spaces, Rep. Math. Phys, 57(2006), 357-366.
\bibitem[13]{13}P. Meessen, Homogeneous Lorentzian spaces admitting a
homogeneous structure of type
 $\tau_{1}\oplus \tau_{3}$, J. Geom.
Phys. 56(2006), 754-761.
\bibitem[14]{14}P. Meessen, Homogeneous Lorentzian spaces whose
null-geodesics are canonically homogeneous, Lett. Math. Phys.
75(2006), 209-212.
\bibitem[15]{15}S. Philip, Penrose limits of homogeneous spaces, J.
Geom. Phys. 56(2006), 1516-1533.


\bibitem[16]{16}G. Randers, On an asymmetrical metric in the
four-space of general relativity, Phys. Rev. 59(1941)195-199.
\bibitem[17]{17}Z. Shen, Differential Geometry of Sprays and Finsler
Space, Kluwer Academic Publishers, 2001.
\bibitem[18]{18}Z. Shen, Volume comparison and its application in
Riemann-Finsler geometry, Advances in Math. 128(1997), 306-328.
\bibitem[19]{19}A. Spiro, Chern`s orthonormal frame bundle of a
Finsler space, Houston. J. Math. 25(1999), 641-659.
\bibitem[20]{20}J. Szenthe, Existence of stationary geodesics of left-invariant
Lagrangians, J. Phys. A:Math. Gen., 34(2001), 165-175.

\end{thebibliography}
\end{document}